\documentclass[12pt]{article}

\usepackage{hyperref}
\usepackage{amsmath, amssymb}
\usepackage{amsfonts, mathrsfs}
\usepackage[cp1251]{inputenc}
\usepackage[active]{srcltx}
\usepackage[ukrainian, russian, english]{babel}

\usepackage{amsfonts, amssymb, mathrsfs}

\begin{document}

\selectlanguage{ukrainian} 

\setcounter{page}{1}

\noindent {\small УДК 517.5} \vskip 3mm

\noindent \textbf{A.S. Serdyuk, I.V. Sokolenko} \\ {\small (Institute of Mathematics
of The National Academy of Sciences of Ukraine, Kiev)}\vskip 3mm

\noindent \textbf{А.\,С. Сердюк, І.\,В. Соколенко} \\  {\small (Iнститут математики Національної академії наук України, Київ)} \vskip 3mm

 \noindent  \textbf{\Large Approximation by linear methods of\\ classes of \boldmath{$(\psi,\bar\beta)-$}differentiable functions}
 \vskip 5mm 

\noindent  \textbf{\Large Наближення лінійними методами\\ класів   \boldmath{$(\psi,\bar\beta)-$}диференційовних функцій}
 \vskip 5mm
 
\noindent  \textbf{\Large Приближение линейными методами\\ классов \boldmath{$(\psi,\bar\beta)-$}дифференцируемых функций}
 \vskip 5mm 
 
{\small \it \noindent We calculate  the least upper bounds for
approximations in the   metric of the space $L_2$ by linear methods of summation of Fourier series   on classes of periodic functions $L^\psi_{\bar\beta,1}$ defined by sequences of multipliers $\psi=\psi(k)$ and shifts of argument  $\bar\beta=\beta_k$.
 \hfill}
 \vskip 3mm
 
{\small \it \noindent Обчислено  точні верхні межі
наближень в метриці простору $L_2$ лінійними методами підсумовування рядів Фур'є   класів періодичних функцій
$L^\psi_{\bar\beta,1}$, які задаються послідовностями мультиплікаторів $\psi=\psi(k)$ та зсувів аргументу $\bar\beta=\beta_k$. \hfill}
\vskip 3mm

{\small \it \noindent Вычислены точные верхние грани приближений в метрике пространства  $L_2$ линейными методами суммирования рядов Фурье классов периодических функций $L^\psi_{\bar\beta,1}$, задающихся последовательностями мультипликаторов $\psi=\psi(k)$ и сдвигов аргумента $\bar\beta=\beta_k$. \hfill}
\vskip 3mm

\normalsize\vskip 3mm

Нехай
$\ C\ $ i $\ L_p, \ 1\le p\le\infty,$ --- простори
$2\pi$-періодичних   функцій  зі стандартними нормами $\ \|\cdot\|_{C}\ $ i  $\ \|\cdot\|_{L_p}$. Одиничну кулю в просторі $L_p$ позначатимемо через $B_p,$ тобто
$$
  B_p=\{\varphi\in L_p:\ \|\varphi\|_{L_p}\le1\},\ \ \ 1\le p\le\infty.
$$

Нехай, далі, $\psi=\psi(k)$ i $\bar\beta=\beta_k$, $k=1,2,\ldots,$ --- довільні послідовності
дійсних чисел такі, що ряд
\begin{equation*}\label{3}
  \sum\limits_{k=1}^\infty \psi(k)\cos\left(kt-\frac{\beta_k\pi}2\right),
\end{equation*}
є рядом  Фур'є деякої функції $\Psi_{\bar\beta}\in L_1.\  $ Тоді через  $L^\psi_{\bar\beta,p}, 1{\le }p{\le}\infty,$  позначимо множину всіх $2\pi$-періодичних функцій $f$, які майже скрізь зображуються за допомогою  згортки
\begin{equation}\label{7'}
    f(x)=\frac{a_0 }{2}+\frac{1}{\pi}\int\limits_{-\pi}^{\pi}
 \varphi(x-t) \Psi_{\bar\beta}(t)dt=\frac{a_0 }{2}+\big(\varphi\ast\Psi_{\bar\beta}\big)(x), \ \ \  a_0\in\mathbb R,  \ \ \  \varphi\in B_p^0.
\end{equation}
де $
  B_p^0=\{\varphi\in B_p:\  \varphi\perp1\}.
$
 Функцію $\varphi$ в зображенні (\ref{7'}) називають $(\psi,\bar{\beta})$-похідною функції $f$ і позначають  $f^{ \psi}_{\bar{\beta}}(x)$. Поняття $(\psi,\bar{\beta})$-похідної введено О.І.~Степанцем
 (див., наприклад, [\ref{Stepanets2002_1}]).  Якщо $\varphi\in B_1^0,\ $ a $\ \Psi_{\bar\beta}\in L_p,$ то з нерівності Юнга для згорток (див. [\ref{Stepanets2002_1}, c. 293]):
$$
\|y\ast z\|_{L_p}\le\frac1\pi\|y\|_{L_s}\|z\|_{L_q},  1\le s\le p\le\infty,\ \frac1q=1-\frac1s+\frac1p, y\in L_s, z\in L_q,
$$
випливає, що $L^\psi_{\bar\beta,1}\subset L_p, 1\le p\le \infty.$ При $p=2$
 включення $\Psi_{\bar\beta}\in L_2$ еквівалентне виконанню умови
\begin{equation}\label{6}
    \sum_{k=1}^\infty \psi^2(k)<\infty.
\end{equation}

 Розглянемо послідовності функцій $\lambda_k(\delta)$ і $\mu_k(\delta)$, які задані на деякій множині    $E\subset \mathbb R$ з граничною точкою $\delta_0$ і задовольняють умови
\begin{equation}\label{7''}
  \begin{array}{lll}
    \lambda_0(\delta)=1,& \mu_0(\delta)=0 & \forall\delta\in E, \\ \\
    \lim\limits_{\delta\rightarrow \delta_0}\lambda_k(\delta)=1,& \lim\limits_{\delta\rightarrow \delta_0}\mu_k(\delta)=0 & k=1,2,\ldots
  \end{array}
\end{equation}
При довільному фіксованому $\delta\in E$ означимо лінійний оператор $U_\delta=U_\delta(\lambda;\mu)$, який кожній функції $f\in L^\psi_{\bar\beta,1}$ ставить у відповідність функцію
\begin{equation*}\label{7}
    U_\delta(\lambda;\mu;f;x)=\frac{a_0}{2}+\sum\limits_{k=1}^\infty(\lambda_k(\delta)(a_k\cos kx+b_k\sin
 kx)+
$$
$$+\mu_k(\delta)(-b_k\cos kx+a_k\sin
 kx)),
\end{equation*}
де $a_k, b_k$ --- коефіцієнти Фур'є функції $f$.

У даній роботі розглядається задача про знаходження точних значень величин
 \begin{equation}\label{9''}
    {\cal E}_\delta(L^\psi_{\bar\beta,1}; \lambda; \mu)_{L_2}=\sup\limits_{f\in L^\psi_{\bar\beta,1}}\|f(x)-U_\delta(\lambda;\mu;f;x)\|_{L_2}.
 \end{equation}
 Має місце наступне твердження.

 \textbf{Теорема 1.}\textit{ Нехай послідовність $\psi=\psi(k)$ задовольняє умову $(\ref{6})$,
  а послідовності функцій $\lambda_k(\delta)$ і $\mu_k(\delta)$ --- умови $(\ref{7''})$. Тоді для довільної послідовності дійсних чисел $\bar\beta=\beta_k $ і довільного $\delta\in E\subset \mathbb R$
 \begin{equation*}\label{2t2}
 {\cal E}_\delta(L^\psi_{\bar\beta,1}; \lambda; \mu)_{L_2}=\frac1{\sqrt\pi}\left(\sum\limits_{k=1}^\infty\left((1-\lambda_k(\delta))^2+\mu^2_k(\delta)\right)\psi^2(k)\right)^{1/2}.
 \end{equation*}
 }

 \textit{\textbf{Доведення.}}  Оскільки згідно з (\ref{7'}) для будь-якої $  f\in L^\psi_{\bar\beta,1}$
\begin{equation*}\label{7_1}
   a_k\cos kx+b_k\sin kx=\frac1\pi\int\limits_{-\pi}^\pi \psi(k)\cos\left(k(x-t)-\frac{\beta_k\pi}2\right)f^{\psi}_{\bar\beta}(t)dt,
\end{equation*}
\begin{equation*}\label{7_1'}
   -b_k\cos kx+a_k\sin
 kx=\frac1\pi\int\limits_{-\pi}^\pi \psi(k)\sin\left(k(x-t)-\frac{\beta_k\pi}2\right)f^{\psi}_{\bar\beta}(t)dt,
\end{equation*}
то майже скрізь виконується рівність
\begin{equation}\label{27_2}
    U_\delta(\lambda;\mu;f;x)=\frac{a_0}{2}+\frac1\pi\int\limits_{-\pi}^\pi f^{\psi}_{\bar\beta}(x-t)U_\delta(\lambda;\mu;t)dt, 
\end{equation}
в якій $U_\delta(\lambda;\mu;t)$ --- функція, ряд Фур'є якої можна представити у вигляді
\begin{equation*}\label{7_3}
     \sum\limits_{k=1}^\infty\psi(k)\left(\lambda_k(\delta)\cos\left(kt{-}\frac{\beta_k\pi}2\right){+}\mu_k(\delta)\sin\left(kt{-}\frac{\beta_k\pi}2\right)\right).
\end{equation*}

  З рівностей (\ref{7'}) i (\ref{27_2})   одержуємо
 $$
   {\cal E}_\delta(L^\psi_{\bar\beta,1}; \lambda; \mu)_{L_2}=
     \sup\limits_{\varphi\in B_1^0}\Bigg\|\frac1\pi\int\limits_{-\pi}^\pi \varphi(x-t)\big(\Psi_{\bar\beta}(t)-U_\delta(\lambda;\mu;t)\big)dt\Bigg\|_{L_2}.
   $$
 Оскільки  функція $\Psi_{\bar\beta}(t)-U_\delta(\lambda;\mu;t)$ ортогональна до будь-якої константи, то
   $$
   {\cal E}_\delta(L^\psi_{\bar\beta,1}; \lambda; \mu)_{L_2}=
  \sup\limits_{\varphi\in B_1}\Bigg\|\frac1\pi\int\limits_{-\pi}^\pi \varphi(x-t)\big(\Psi_{\bar\beta}(t)-U_\delta(\lambda;\mu;t)\big)dt\Bigg\|_{L_2}.
   $$

  Як відомо (див., наприклад, [\ref{Kornejchuk1987}, наслідок Д1.2., c. 392]), якщо $u\in L_p,\ 1\le p\le\infty,$ то
 \begin{equation}\label{12'''}
 \|u \|_{L_p}=\sup_{g\in B_{p'}} \int\limits_{-\pi}^\pi g(t)u(t)dt, \ \ \ \frac1p+\frac1{p'}=1.
 \end{equation}

  Враховуючи    співвідношення   (\ref{12'''}), інваріантність множини $B_2$ відносно зсуву аргументу та рівність Парсеваля, маємо
   $$
   {\cal E}_\delta(L^\psi_{\bar\beta,1}; \lambda; \mu)_{L_2}=
   \frac1\pi\sup\limits_{\varphi\in B_1}\sup_{g\in B_2} \int\limits_{-\pi}^\pi g(x) \int\limits_{-\pi}^\pi \varphi(t)\big(\Psi_{\bar\beta}(x-t)-U_\delta(\lambda;\mu;x-t)\big)dtdx =
   $$
   $$
   =\frac1\pi\sup\limits_{g\in B_2}\sup_{\varphi\in B_1 } \int\limits_{-\pi}^\pi \varphi(t) \int\limits_{-\pi}^\pi g(x+t)\big(\Psi_{\bar\beta}(x)-U_\delta(\lambda;\mu;x)\big)dxdt =
   $$
   $$
   =\frac1\pi\sup\limits_{g\in B_2}
   \Bigg\|  \int\limits_{-\pi}^\pi g(x+t)\big(\Psi_{\bar\beta}(x)-U_\delta(\lambda;\mu;x)\big)dx\Bigg\|_{L_\infty} =
   $$
   $$
   =\frac1\pi\sup\limits_{g\in B_2}\int\limits_{-\pi}^\pi g(x)\big(\Psi_{\bar\beta}(x)-U_\delta(\lambda;\mu;x)\big)dx=
   \frac1\pi\big\|\Psi_{\bar\beta} -U_\delta(\lambda;\mu)\big\|_{L_2}=
   $$
   $$
   {=}\frac1\pi\bigg\|
   \sum\limits_{k=1}^\infty\psi(k)\!\left(\!(1{-}\lambda_k(\delta))\cos\left(kx{-}\frac{\beta_k\pi}2\right){+}
  \mu_k(\delta)\sin\left(kx{-}\frac{\beta_k\pi}2\right)\!\right)\!\bigg\|_{L_2}\! {=}
   $$
  $$
=\frac1{\sqrt\pi}\left(\sum\limits_{k=1}^\infty\left((1-\lambda_k(\delta))^2+\mu^2_k(\delta)\right)\psi^2(k)\right)^{1/2}.
$$

 Теорему  доведено.

\textit{\textbf{Зауваження 1.}} З теореми 1 даної роботи і теореми 1 з роботи [\ref{Serdyuk_Sokolenko_2011}]   випливає, що
 $$
    {\cal E}_\delta(L^\psi_{\bar\beta,1}; \lambda; \mu)_{L_2}={\cal E}_\delta(C^\psi_{\bar\beta,2}; \lambda; \mu)_C,
 $$
де
$
{\cal E}_\delta(C^\psi_{\bar\beta,2}; \lambda; \mu)_C=\sup\limits_{f\in C^\psi_{\bar\beta,2}}\|f(x)-U_\delta(\lambda;\mu;f;x)\|_C,
$
 а $C^\psi_{\bar\beta,2}=C\cap L^\psi_{\bar\beta,2}.$

Таким чином,  мають місце всі твердження по  наближенню класів $L^\psi_{\bar\beta,1}$ в метриці простору $L_2, $ які є аналогічними до тверджень, що стосуються рівномірних наближень класів $C^\psi_{\bar\beta,2}$, отриманих у  [\ref{Serdyuk_Sokolenko_2011}].
Наведемо деякі з них.

Розглянемо лінійні  поліноміальні методи наближення функцій з класів $L^\psi_{\bar\beta,1}$. Нехай $\Lambda=\|\lambda_k^{(n)}\|\ $ i $\ {\rm M}=\|\mu_k^{(n)}\|,$  $n=0,1,\ldots,  $ \mbox{$\ k=0,1,\ldots,$~---} нескінченні трикутні матриці чисел такі, що:
\begin{equation}\label{a1}
   \begin{array}{lll}
    \lambda_0^{(n)}=1,& \mu_0^{(n)}=0, &  n=0,1,2,\ldots,\\ \\
    \lambda_k^{(n)}=0, &  \mu_k^{(n)}=0, & k=n+1, n+2, \ldots, \\ \\
    \lim\limits_{n\rightarrow \infty}\lambda_k^{(n)}=1,& \lim\limits_{n\rightarrow \infty}\mu_k^{(n)}=0, & k=1,2,\ldots,
  \end{array}
\end{equation}
і лінійний оператор $U_n=U_n(\Lambda;{\rm M})$ кожній функції $f\in L_1$
ставить у відповідність тригонометричний поліном
\begin{equation}\label{a2}
    U_n(\Lambda;{\rm M};f;x)=\frac{a_0}{2}+\sum\limits_{k=1}^n(\lambda_k^{(n)}(a_k\cos kx+b_k\sin
 kx)+
$$
$$+\mu_k^{(n)}(-b_k\cos kx+a_k\sin kx)),
\end{equation}
де $a_k\ i\ b_k$ --- коефіцієнти Фур'є функції $f.$

 Для  методів $U_n$ теорема 1 формулюється таким чином.

 \textbf{Теорема \boldmath{$1'$.}}\textit{ Нехай послідовність $\psi(k)$ задовольняє умову $(\ref{6})$, а  $\Lambda=\|\lambda_k^{(n)}\|\ $ i $\ {\rm M}=\|\mu_k^{(n)}\|$ ---  умови $(\ref{a1})$. Тоді для довільної послідовності дійсних чисел $\bar\beta=\beta_k$ і довільного $n\in \mathbb N$
  \begin{equation}\label{1t'}
 {\cal E}_n(L^\psi_{\bar\beta,1}; \Lambda;{\rm M})_{L_2}=\sup\limits_{f\in L^\psi_{\bar\beta,1}}\|f(x)-U_n(\Lambda;{\rm M};f;x)\|_{L_2}=
$$
$$=\frac1{\sqrt\pi}\left(\sum\limits_{k=1}^n\left((1-\lambda_k^{(n)})^2+(\mu_k^{(n)})^2\right)\psi^2(k)+\sum\limits_{k=n+1}^\infty\psi^2(k)\right)^{1/2}.
 \end{equation}
 }

З теореми $1'$  для деяких класичних лінійних методів наближення випливають такі наслідки.

 Якщо
\begin{equation*}\label{n1}
    \begin{array}{ccc}
      \lambda_k^{(n)}=\left\{
                   \begin{array}{ll}
                     1, & 0\le k\le n, \\
                     0, & k>n,
                   \end{array}
                 \right. & i & \mu_k^{(n)}\equiv0,
    \end{array}
\end{equation*}
то тригонометричний поліном $U_n(\Lambda;{\rm M};f;x)$ є частинною сумою
Фур'є  $S_n(f;x)$ функції $f$ порядку $n$. В цьому випадку
 \begin{equation}\label{1n1}
 {\cal E}_n(L^\psi_{\bar\beta,1}; \Lambda;{\rm M})_{L_2}={\cal E}(L^\psi_{\bar\beta,1}; S_n)_{L_2}
=\frac1{\sqrt\pi}\left(\sum\limits_{k=n+1}^\infty\psi^2(k)\right)^{1/2}.
 \end{equation}

Для класів $L^\psi_{\bar\beta,1}$   при $\psi(k)=q^{k},$  $0< q<1,$ з  (\ref{1n1}) одержуємо
\begin{equation}\label{2n1''}
{\cal E}(L^\psi_{\bar\beta,1}; S_n)_{L_2}=\frac{q^{n+1}}{\sqrt{\pi(1-q^2)}}.
 \end{equation}
При $\beta_k\equiv \beta, \  \beta\in\mathbb R,$ рівність (\ref{2n1''}) уточнює асимптотичну рівність (12) з [\ref{Serdyuk2005_10}] у тому сенсі, що зазначена рівність  (12) з [\ref{Serdyuk2005_10}] при $p=2$ залишиться вірною, якщо в ній вилучити залишковий член.

З   (\ref{1t'}) i (\ref{1n1}) випливає, що найкращим лінійним  методом наближення вигляду (\ref{a2})  класів $L^\psi_{\bar\beta,1}$ у метриці простору $L_2$ є метод Фур'є. Цей факт можна отримати з інших міркувань. Нехай $f\in L_p,$ $ 1\le p\le\infty,$ i
\begin{equation*}\label{a3}
    E_n(f)_{L_p}=\inf\limits_{T_n}\|f-T_n\|_{L_p}
\end{equation*}
--- найкраще наближення функції $f$ в метриці простору $L_p$ тригонометричними поліномами, порядку не вищого ніж $n$.

Нехай, далі, $\mathfrak N \subset L_p$ і $U_n$ --- лінійний оператор вигляду (\ref{a2}). 
Тоді величину
\begin{equation}\label{a4}
    {\cal E}_n(\mathfrak N)_{L_p}=\inf\limits_{U_n}\sup\limits_{f\in\mathfrak N }\|f-U_n(f)\|_{L_p}
\end{equation}
називають найкращим лінійним  наближенням класу $\mathfrak N$ за допомогою лінійних операторів вигляду (\ref{a2})
у метриці простору $L_p,$ а оператор $U_n^\ast$, який реалізує $\inf$ у правій частині (\ref{a4}), називається найкращим лінійним оператором наближення класу $\mathfrak N.$

Виконується наступне твердження.

\textbf{Теорема 2.}\textit{ Нехай $1\le p\le \infty$, а послідовності $\psi=\psi(k)\ i \ \bar\beta=\beta_k$ такі, що $\psi(k)\neq0\ \forall k\in\mathbb N\ $ i $\ \Psi_{\bar\beta}\in L_{p}$. Тоді для довільного $n\in\mathbb N$
  \begin{equation}\label{2t}
 {\cal E}_n(L^\psi_{\bar\beta,1})_{L_p}=\frac1{\pi}E_n(\Psi_{\bar\beta})_{L_p}.
 \end{equation}
 }

  \textit{\textbf{Доведення.}} Неважко переконатися, що  якщо $\psi(k){\neq}0 \ \forall k{\in}\mathbb N$, то для будь-якої $  f\in L^\psi_{\bar\beta,1}$  і довільного тригонометричного полінома з нульовим середнім значенням  $T_n^0(t)=\sum\limits_{k=1}^n(\alpha_k\cos kt+\gamma_k\sin kt),$  $\alpha_k,\ \gamma_k\in\mathbb R,$  згортка
\begin{equation}\label{b1}
     \frac{a_0}{2}+\int\limits_{-\pi}^\pi f^\psi_{\bar\beta}(x-t)T_n^0(t)dt
\end{equation}
 зображується у вигляді полінома $U_n(\Lambda;{\rm M};f;x)$, означеного в (\ref{a2}), при
\begin{equation*}\label{b3}
    \lambda_k^{(n)}=\frac{\alpha_k}{\psi(k)}\cos\frac{\beta_k\pi}2+\frac{\gamma_k}{\psi(k)}\sin\frac{\beta_k\pi}2,
\end{equation*}
\begin{equation*}\label{b4}
    \mu_k^{(n)}=\frac{\gamma_k}{\psi(k)}\cos\frac{\beta_k\pi}2-\frac{\alpha_k}{\psi(k)}\sin\frac{\beta_k\pi}2;
\end{equation*}
і навпаки, довільний тригонометричний поліном $U_n(\Lambda;{\rm M};f;x)$ вигляду (\ref{a2}) зображується у  вигляді згортки (\ref{b1}) з поліномом $T_n^0(t),$ коефіцієнти якого мають вигляд
\begin{equation*}\label{b3'}
    \alpha_k={\psi(k)}\left(\lambda_k^{(n)}\cos\frac{\beta_k\pi}2-\mu_k^{(n)}\sin\frac{\beta_k\pi}2\right),
\end{equation*}
\begin{equation*}\label{b4'}
    \gamma_k={\psi(k)}\left(\lambda_k^{(n)}\sin\frac{\beta_k\pi}2+\mu_k^{(n)}\cos\frac{\beta_k\pi}2\right).
\end{equation*}

Тому, з урахуванням (\ref{7'}), (\ref{a2}) i (\ref{a4}), маємо
 $$
  {\cal E}_n(L^\psi_{\bar\beta,1})_{L_p}=\inf\limits_{U_n}\sup\limits_{f\in L^\psi_{\bar\beta,1}}\|f-U_n(f)\|_{L_p}=
$$
$$
=
\inf\limits_{T_n^0}\sup\limits_{f\in L^\psi_{\bar\beta,1}}
    \left\|\frac{1}{\pi}\int\limits_{-\pi}^{\pi}
 f^\psi_{\bar\beta}(x-t) \left(\Psi_{\bar\beta}(t)-T_n^0(t)\right)dt\right\|_{L_p}=
 $$
$$
=
\inf\limits_{T_n^0}\sup\limits_{\varphi\in B_1^0}
    \left\|\frac{1}{\pi}\int\limits_{-\pi}^{\pi}
 \varphi(x-t) \left(\Psi_{\bar\beta}(t)-T_n^0(t)\right)dt\right\|_{L_p}.
 $$
Враховуючи, що $(\Psi_{\bar\beta} -T_n^0)\perp1$, співвідношення  (\ref{12'''})  та інваріантність  $B_{p'}$ відносно зсуву аргумента, одержуємо
 $$
{\cal E}_n(L^\psi_{\bar\beta,1})_{L_p}=
\frac{1}{\pi}\inf\limits_{T_n^0}\sup\limits_{\varphi\in B_1}\sup\limits_{g\in B_{p'}}
      \int\limits_{-\pi}^{\pi}g(x)\int\limits_{-\pi}^{\pi}
 \varphi(t) \left(\Psi_{\bar\beta}(x-t)-T_n^0(x-t)\right)dtdx=
   $$
$$
=\frac{1}{\pi}\inf\limits_{T_n^0}\sup\limits_{g\in B_{p'}}\sup\limits_{\varphi\in B_1}
      \int\limits_{-\pi}^{\pi} \varphi(t)\int\limits_{-\pi}^{\pi}
g(x) \left(\Psi_{\bar\beta}(x-t)-T_n^0(x-t)\right)dxdt=
$$
$$
=\frac{1}{\pi}\inf\limits_{T_n^0}\sup\limits_{g\in B_{p'}}\bigg\| \int\limits_{-\pi}^{\pi}
g(x+t) \left(\Psi_{\bar\beta}(x)-T_n^0(x)\right)dx\bigg\|_{L_\infty}=
$$
$$
=\frac{1}{\pi}\inf\limits_{T_n^0}\sup\limits_{g\in B_{p'}} \int\limits_{-\pi}^{\pi}
g(x) \left(\Psi_{\bar\beta}(x)-T_n^0(x)\right)dx.
$$

Як відомо з [\ref{Kornejchuk1987}, c. 27], для довільної функції $\ u\in L_{p}, 1\le p\le\infty,$ має місце співвідношення двоїстості:
 \begin{equation}\label{12''}
 \inf_{\alpha\in \mathbb R}\|u-\alpha \|_{L_p}=\sup_{g\in B_{p'}^0} \int\limits_{-\pi}^\pi g(t)u(t)dt, \ \ \ \frac1p+\frac1{p'}=1.
 \end{equation}

Оскільки $(\Psi_{\bar\beta} -T_n^0)\perp1$, то в силу   (\ref{12''})
$$
   {\cal E}_n(L^\psi_{\bar\beta,1})_{L_p}= \frac{1}{\pi}\inf\limits_{T_n^0}\sup\limits_{g\in B_{p'}^0}
      \int\limits_{-\pi}^{\pi}
 g(x) \left(\Psi_{\bar\beta}(x)-T_n^0(x)\right)dx =
 $$
 $$
 =\frac{1}{\pi}\inf\limits_{T_n^0}\inf\limits_{\alpha\in\mathbb R}\|\Psi_{\bar\beta}(t)-T_n^0(t)-\alpha\|_{L_p}=\frac{1}{\pi}E_n(\Psi_{\bar\beta})_{L_p}.
$$
 Теорему 2 доведено.

 З міркувань, використаних при доведенні теореми 2, випливає також, що при $\psi(k)\neq0$ $\ \forall k\in\mathbb N$, для класу $L^\psi_{\bar\beta,1}$ серед усіх лінійних методів наближення  $U_n$ вигляду (\ref{a2}) у метриці простору $L_p$ найкращим є метод $U_n^\ast(\Lambda;{\rm M};f;x)$, що породжений системою чисел
\begin{equation*}\label{b3}
    \lambda_k^{(n)}=\frac{\alpha_k^\ast}{\psi(k)}\cos\frac{\beta_k\pi}2+\frac{\gamma_k^\ast}{\psi(k)}\sin\frac{\beta_k\pi}2,
\end{equation*}
\begin{equation*}\label{b4}
    \mu_k^{(n)}=\frac{\gamma_k^\ast}{\psi(k)}\cos\frac{\beta_k\pi}2-\frac{\alpha_k^\ast}{\psi(k)}\sin\frac{\beta_k\pi}2,
\end{equation*}
де $\alpha_k^\ast$ i $\gamma_k^\ast$ --- коефіцієнти полінома найкращого наближення твірного ядра $\Psi_{\bar\beta}$ у метриці простору $L_{p}.$
Оскільки $E_n(\Psi_{\bar\beta})_2=\|\Psi_{\bar\beta}- S_n(\Psi_{\bar\beta})\|_2,$ то найкращим лінійним  методом наближення $U_n^\ast$  класів $L^\psi_{\bar\beta,1}$ у метриці простору $L_2$ є метод Фур'є.

Наведемо  наслідок з теореми $1'$ для методу Валле Пуссена.
Нехай $n=0,1,2,\ldots, $ $0\le m\le n,\ m\in\mathbb N$ i
 \begin{equation*}\label{n2}
\begin{array}{cc}
      \lambda_k^{(n)}=\left\{
                   \begin{array}{ll}
                     1, & 0\le k\le n-m, \\
                     1-\frac{k-n+m}{m+1}, & n-m<k\le n,\\
                     0, & k>n,
                   \end{array}
                 \right. & \mu_k^{(n)}\equiv0.
    \end{array}
\end{equation*}
Тоді поліном $U_n(\Lambda;{\rm M};f;x)$ вигляду (\ref{a2}) є сумою Валле Пуссена  $V_{n,m}(f;x)$ функції $f$ i згідно з (\ref{1t'})
 \begin{equation}\label{1n2}
 {\cal E}(L^\psi_{\bar\beta,1}; \Lambda;{\rm M})_{L_2}={\cal E}(L^\psi_{\bar\beta,1};V_{n,m})_{L_2}=
$$
$$
=\frac1{\sqrt\pi} \left(\frac1{(m{+}1)^2}\sum\limits_{k{=}n{-}m{+}1}^{n}(k{-}n{+}m)^2\psi^2(k)+
\sum\limits_{k=n+1}^\infty\psi^2(k)\right)^{1/2},
 \end{equation}

В [\ref{Serdyuk2010}, с. 1680] доведено, що  при   $\ 0<q<1$
\begin{equation}\label{2n1'}
    \frac1{(m+1)^2}\sum\limits_{k{=}n{-}m{+}1}^{n}(k{-}n{+}m)^2 q^{2k}+\sum\limits_{k=n+1}^\infty q^{2k}=
$$
$$
=\frac{q^{2(n-m+1)}(1+q^2- q^{2(m+1)}(2m+3-q^2(2m+1)))}{(m+1)^2(1-q^2)^3}.
\end{equation}
З (\ref{1n2}) i (\ref{2n1'}) для класів   $L^q_{\bar\beta,1}$  отримуємо рівність
\begin{equation}\label{101}
{\cal E}(L^q_{\bar\beta,1};V_{n,m}
)_{L_2}=
\frac{q^{n-m+1}}{\sqrt\pi (m+1)}\sqrt{\frac{1+q^2-q^{2(m+1)}(2m+3-q^2(2m+1))}{(1-q^2)^3}}.
\end{equation}

При $\beta_k\equiv \beta, \  \beta\in\mathbb R,$ рівність (\ref{101}) уточнює асимптотичну рівність (71), отриману в [\ref{Serdyuk2010}]. А саме зазначена  рівність (71) з [\ref{Serdyuk2010}] залишається вірною, якщо в ній вилучити  залишковий член.

\footnotesize
\begin{enumerate}
\item\label{Stepanets2002_1} {\it Степанец А.И.\/} Методы теории приближений: В
2 ч. --- Киев: Ин-т математики НАН Украины, 2002. --- Ч.\,1. ---
427\,c.

\item \label{Kornejchuk1987} {\it Корнейчук Н.П.\/} Точные константы
в теории приближения. --- М.: Наука, 1987. --- 423~с.

\item  \label{Serdyuk_Sokolenko_2011}  {\it Сердюк А.С., Соколенко І.В.\/} Рівномірні наближення класів $(\psi,\bar\beta)$-диференційовних функцій лінійними методами // Проблеми теорії наближення функцій та суміжні питання: Зб. праць Інституту математики  НАН  України. --- 2011. ---   Т. 8, №1.--- С.~181-189.

\item \label{Serdyuk2005_10} {\it  Сердюк А.С.\/} Наближення класів аналітичних функцій сумами Фур’є в метриці простору $L_p$. //
Укр. мат. журн. --- 2005. --- {\bf 57}, № 10. --- C.1395–1408.

\item \label{Serdyuk2010} {\it  Сердюк А.С.\/}
Приближение интегралов Пуассона суммами Валле Пуссена в равномерной и интегральной метриках //
Укр. мат. журн. --- 2010. --- {\bf 62}, № 12. --- C.~1672–1686.

\end{enumerate}
\label{end}

\textbf{Contact information:}

Department of the Theory of Functions, Institute of Mathematics
of The National Academy of Sciences of Ukraine, 3, Tereshenkivska st., 01601, Kyiv, Ukraine.

E-mail: sokol@imath.kiev.ua, serdyuk@imath.kiev.ua

 \end{document}